# Design and control of cross-coupled mechanisms driven by AC brushless servo motors

AM CONNOR
University of Bath, UK

**Abstract**
*This paper presents an overview of a design methodology for the optimal synthesis of hybrid mechanisms. Hybrid mechanisms have been defined as multi-degree of freedom systems where the input motions are supplied by different motor types. In this work a five bar mechanism is designed for a given task under the constraint that one input axis rotates with constant velocity whilst the other input can exhibit any motion requirement. A machine of this type is classified as being cross-coupled due to the mechanical linkage between the input axes. Cross-coupling implies that the input motion on one axis effects the position of the other input axis. This can lead to either opposition to, or accentuation of the control system input. Such a system as this is difficult to control due to the compensation for this on each axis leading to further disturbance. Results are presented for a real machine operating in this way and the actual output of the machine is compared to the desired input of the machine.*

## NOTATION

| | | | |
|---|---|---|---|
| p,q,r,s | Mechanism link lengths | $\theta_j$ | Mechanism angles |
| k | Number of precision points | $K_P$ | Proportional gain |
| $x_d$ | Desired x coordinate of end effector | $K_I$ | Integral gain |
| $x_a$ | Actual x coordinate of end effector | $K_D$ | Differential gain |
| $y_d$ | Desired y coordinate of end effector | $K_V$ | Velocity feedback gain |
| $y_a$ | Actual y coordinate of end effector | $K_F$ | Velocity feedforward gain |
| $a_j, b_j$ | Fourier coefficents | $e_i$ | Position error |
| n | Number of harmonic coefficients | $d_i$ | Demand position |
| m | Number of immobile closures | $p_i$ | Measured position |

# 1 INTRODUCTION

The nature of the modern production environment has introduced the demand that flexibility must be high and that down time must be low for a manufacturer to remain competitive. This has resulted in changing the way in which machinery is designed. Conventional design methods consist of inserting cams and linkages into the driving mechanisms. Whilst a degree of flexibility is given by having interchangeable parts, this results in high down times. However, traditional machines have good transmission and the potential for energy regeneration.

In an attempt to achieve high flexibility and low down time, programmable electric motors have been used to drive the output directly, as opposed to through a linkage. However, for applications that follow non-uniform trajectories and demand high rates of change of torque the motors generate considerable heat and the dissipation of this heat is manifested as a high power rating.

Generally, a programmable motor used in conjunction with closed loop control is referred to as a servo motor. Servo motor systems have high flexibility but may suffer from poor energy regeneration. Constant velocity (CV) motors, and traditional machines, have good transmission and the potential for good energy regeneration. This is often achieved by the inclusion of a flywheel on the drive shaft. An approach that combines the desirable characteristics of both systems is the use of a hybrid machine.

## 1.1 Hybrid Machines

The initial concept of a hybrid machine investigated by Tokuz (1,2), who developed a system in which the outputs of a constant velocity motor and a programmable servo motor were combined through a differential gearbox to drive a slider-crank mechanism. Figure 1.1 shows a schematic diagram of the machine. The aim of the work was to consider the efficiency of the hybrid rig with the case where the servo motor drives the slider-crank directly. This is an example of a programmable machine, as opposed to a hybrid machine.

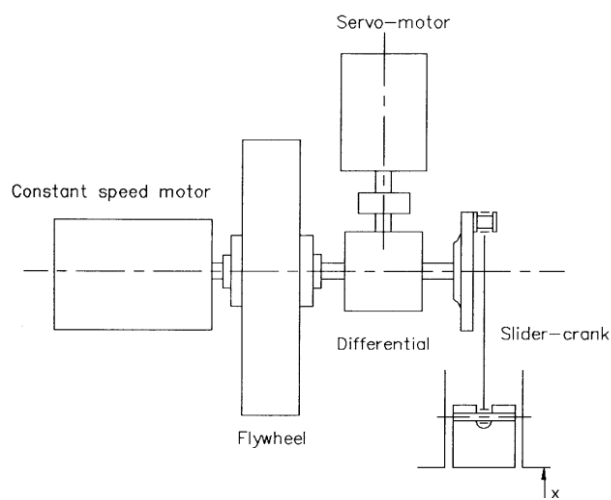

Figure 1.1 : Hybrid arrangement

Several different slider motions were investigated, including a Rise-Return (R-R), a Rise-Dwell-Return (R-D-R) and a Rise-Return-Dwell (R-R-D). For the R-R motion, the hybrid approach showed considerable power savings over the programmable machine. However, for the other motions the hybrid machine required more power to produce the output motion. This can be attributed to the dwell in both the R-D-R and the R-R-D motions. To obtain a dwell, the servo motor must oppose the motion of the constant velocity motor directly which leads to high power consumption.

Later work by Greenough (3) and Bradshaw (4,5) attempted to overcome the difficulties of producing complex output motions by replacing the differential gearbox with a multi-degree of freedom linkage to act as a non-linear gearbox. In this arrangement, the servo motor need not directly oppose the constant velocity motor to obtain a dwell. Their work investigated the case of a hybrid machine designed to have two fully rotational angular inputs and one fully rotational angular output. After researching several alternatives, a Svoboda type II mechanism was used.

The aim of the research was to investigate issues concerning the design, control and applications of hybrid machines of this configuration. One particular application was for a cut to length machine. Figure 1.2 shows such a machine. The machine consists of a pair of high inertia cylinders, each of which is driven by a programmable servo motor. Each cylinder has a blade attached to it. A shearing action is produced when the blades meet which cuts the material. By matching the speed of the blade to the linear speed of the web, a high quality cut is produced.

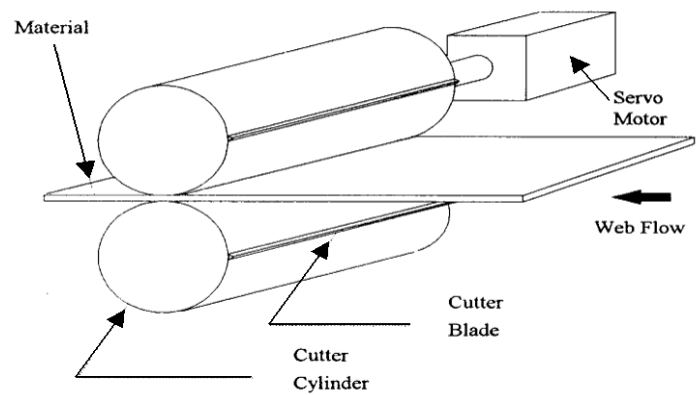

Figure 1.2 : Cut to length machine

Replacing each servo motor with the hybrid machine module reduced the peak power by approximately 70%. There was also a slight increase in speed of operation. The machine module was also tested on a number of intermittent motions, such as the R-D-R and R-R-D motions. In these cases the power savings were less significant.

The work presented in this paper is part of an on going study to investigate further aspects of the hybrid machine concept. The primary aim was to develop a more robust and more efficient design methodology. This will be briefly described in the next section but is fully documented in (6,7).

## 2 HYBRID MACHINE DESIGN METHODOLOGY

The design of mechanisms consists of two stages of synthesis. The first stage is *type* synthesis where the designer selects a mechanism configuration that is likely to produce the desired output motion. The second stage of the process is *dimensional* synthesis where the parameter values for the mechanism link lengths are selected to ensure that the output motion requirement is matched exactly.

This work deals purely with the dimensional synthesis of a hybrid five bar mechanism that is intended to produce a closed curve as an output motion. Figure 2.1 shows the notation in use for the analysis of the five bar mechanism. Each 'stick' in the diagram represents a link between revolute joints. The link t is the virtual ground link, and so all motions in the mechanism are relative to this fixed datum. The link *p* is the CV input and so rotates fully around the ground point. The link *s* is the servo motor input and has a programmable motion.

As the five bar is a two degree of freedom mechanism, it requires two parameters to be defined if the mechanism is to be fully analysed. In the analysis used, these two parameters are $\theta_2$ and $\theta_3$. $\theta_2$ is the input angle associated with the CV motor position and as such is known. $\theta_3$ can be calculated for each position of the CV motor given the desired position of the end effector. This is calculated by defining a vector from the end of link *p* to the desired position. This is illustrated in Figure 2.2.

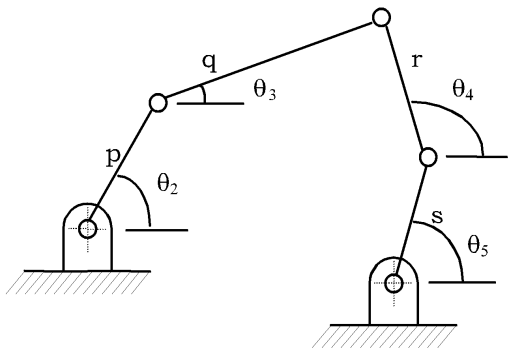
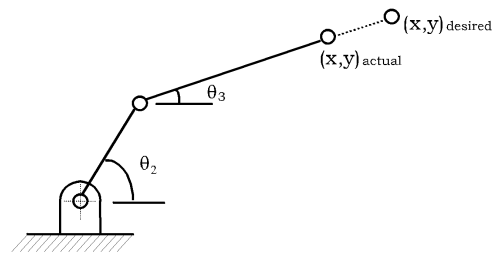

Figure 2.1 : Five bar mechanism                Figure 2.2 : Calculation of $\theta_3$

### 2.2.1 Structural Error
As the length of link *p* is known the actual position of the end effector can be calculated which then allows the structural error to be calculated which is merely the difference between the desired and actual positions. This has difference has been illustrated in diagrammatic form in Figure 2.2. The structural error can be summed around the cycle of the mechanism and then raised to a power to provide a useful merit indicator for use in the automatic dimensional synthesis.

However, there are several other issues that need to be considered before an objective function can be developed. So far only the CV motor input dyad has been considered. Now that values for $\theta_2$ and $\theta_3$ are known it is possible to calculate $\theta_4$ and $\theta_5$. A complication arises that the equations which describe the configuration of the mechanism are quadratic in nature which results in two possible closures for each position with two corresponding values of $\theta_5$.

A closure tracking algorithm has been developed which chooses the best position of link *s* for each position of the CV input. This algorithm is fully described in (6,7,8). Once the complete motion of link *s* is known three other merit terms can be calculated for use in the synthesis.

### 2.2.2 Mobility Penalty

The first of the additional terms used in the synthesis objective function is the mobility penalty function. The actual position of the end effector has already been calculated considering it to be a part of the CV input dyad. Similarly, the remaining two mobile links, *r* and *s*, form another dyad which is used to calculate a penalty function criteria of the overall objective function.

This penalty function is based on the mobility of the mechanism. For a truly mobile mechanism, the dyad formed by the links *r* and *s* should be able to 'close' for all given positions of the CV input crank. This means that the position of the common revolute joint, when considered as part of the dyad formed by *r* and *s*, should be able to reach the actual position defined by the dyad formed by the links *p* and *q*.

### 2.2.3 Harmonic Content and Swept Area

The two remaining terms used in the objective function are the harmonic content of the servo motor motion profile and the swept area of link *s*. These two terms are used to approximate the dynamic performance of the mechanism and force the synthesis method to produce solutions which exhibit desirable a torque distribution without the necessity of carrying out complex dynamic calculations.

The concept of motion swept area was inspired by early work (8), where minimising the RMS of the motor displacements forced the search to locate mechanisms where the links in the closing dyad were very much longer than those in the input dyad. A mechanism with such link length ratios is unlikely to produce even an approximation to good dynamic performance due to the large torque requirements of the servo motor. In calculating the swept area term, the RMS value of the servo motor displacement is multiplied by the length of link *s*. By doing this, an implicit compromise is achieved between the magnitude of the displacements and the length of link *s*.

In addition to investigating the use of a swept area term in the objective function, it was also decided to use a term based on the harmonic content of the motion profile. This decision was also based on previous work (8), where minimising the RMS of the motion produced a motion profile which exhibited considerable oscillation. Penalising profiles with high magnitudes in high order harmonics should direct the search towards a solution where the servo displacement profile is trending towards simple harmonic motion. Calculating the harmonic content of a motion profile can be achieved by utilising a Fourier transform.

### 2.2.4 Synthesis Method

The dimensional synthesis of mechanisms lends itself to the use of numerical optimisation techniques in order to automate the process. The method utilised in this work is a Genetic

Algorithm. Genetic Algorithms (GAs) are a non-derivative based optimisation method based on the Darwinian theory of natural selection and are a broad and effective search. However, GAs have been shown to lack local search power so the synthesis method adopted is a two stage method where the solution found by the GA is passed into a local hill climbing search in order to find the true optimum solution in the region located by the GA. The hill climbing search used was a simple steepest descent algorithm.

### 2.2.4.1 Objective Function

There are four components to the objective function used in this work. The intention is to provide a degree of approximation to a true dynamic objective function, incorporating power consumption, but by using only kinematic terms which are considerably less computationally expensive. The four components are given by;

$$obj_{error} = \left( \sum_{i=1}^{k} \left[ \sqrt{(x_a - x_d)^2 + (y_a - y_d^2)} \right] \right)^2 \quad \ldots(1)$$

$$obj_{harm} = \sum_{i=1}^{n} \left( \sqrt{a_n^2 + b_n^2} \right)^{i+1} \quad \ldots(2)$$

$$obj_{swept} = s \sqrt{\sum_{i=1}^{k} \theta_{5i}^2} \quad \ldots(3)$$

$$obj_{mob} = m^3 \quad \ldots(4)$$

The four components of the objective function are combined into a single numerical value by using a weighted sum approach. The component weightings and the GA control parameters were determined by experimentation and are given below.

| GA Control Parameters | Objective Function Weightings |
|---|---|
| Population size : 40 | Error : 1.0 |
| Crossover rate : 0.85 | Mobility : 1.0 |
| Mutation rate : 0.03 | Swept area : 0.75 |
| | Harmonic content : 0.5 |

The GA utilised standard roulette wheel selection and a linear fitness scaling routine to ensure that the population of solutions did not suffer from premature convergence to a sub-optimal solution.

## 2.3 Results

Using the weightings and control parameters from the previous section a number of sample runs of the method were carried out. In these runs the design parameters were the lengths of links *p*,*q*,*r*, and *s* as well as the (x,y) coordinates of both the CV and the servo motor. The best mechanism located had the following dimensions;

p = 150  q = 250  r = 300  s = 150
CV = (0,0)        Servo = (250,0)

The servo motor input motion displacement profile of this mechanism is shown in Figure 2.3. The velocity and acceleration profiles are shown in Figure 2.4 and Figure 2.5 respectively. The motion of the end effector is shown in Figure 2.6 and the torque requirements of the motors are shown in the graph in Figure 2.7.

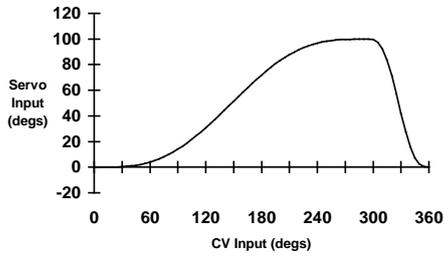

Figure 2.3 : Servo displacement

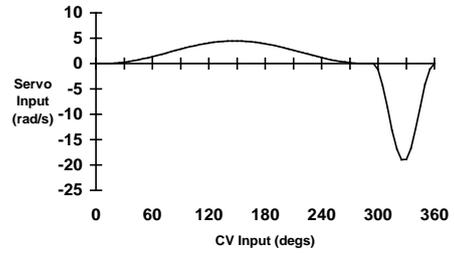

Figure 2.4 : Servo velocity

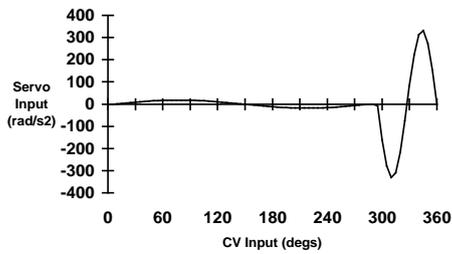

Figure 2.5 : Servo motor acceleration

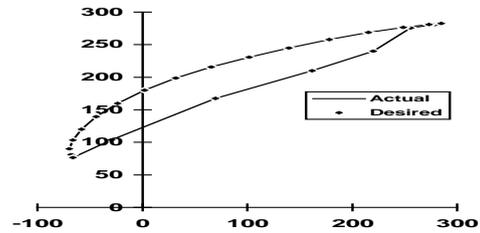

Figure 2.6 : End effector motion

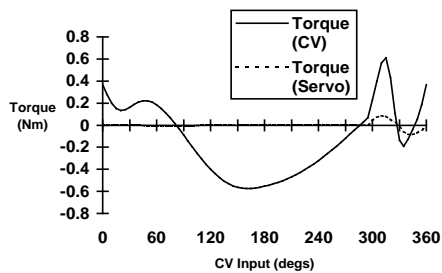

Figure 2.7 : Torque requirements

The minimum, maximum and RMS values for the torque requirement of each motor are given in the table in Table I.

|  | $T_{MIN}$ (Nm) | $T_{MAX}$ (Nm) | $T_{RMS}$ (Nm) |
|---|---|---|---|
| CV Motor | -0.5739 | 0.6107 | 0.6043 |
| Servo Motor | -0.0859 | 0.0862 | 0.0461 |

Table I : Motor torque values

These values indicate that the solution is likely to be acting as a true hybrid device as the torque requirement for the CV motor is very much larger than that for the servo motor. A practical machine was built to verify this conclusion and to compare the actual output motion with the desired output motion.

## 3    PRACTICAL MACHINE

The experimental apparatus consists of two servo motors fixed in position. These two motors are joined together by a five bar linkage of dimensions given in section 2.3. One of these motors is forced to act as a CV motor whilst the other is programmed to follow a given motion profile. The output of the end effector is traced using a camera. An LED is located on the end effector and a long exposure time used to ensure that the trajectory is recorded on film.

The sketch in Figure 3.1 shows the arrangement of the motors, housing, misalignment coupling and bearings. The linkage mechanism is attached to the output shafts. A sketch of the linkage in the starting position is given in Figure 3.2.

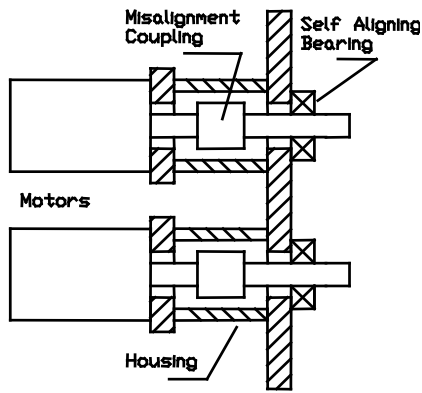
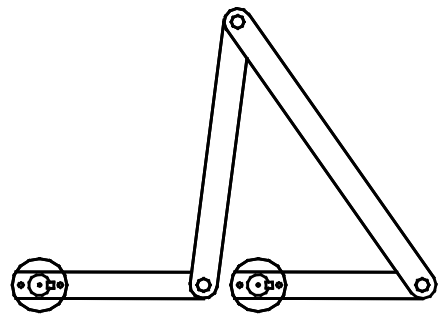

Figure 3.1 : Schematic arrangement     Figure 3.2 : Five bar linkage

Essentially, each servodrive consists of a motor with low rotor inertia and an amplifier which is used to modulate the potential difference across the motor coils so as to produce a variable speed output. A controller is used to interpret feedback information and provide command signals to modulate the potential difference as required. The control system operates by sampling the position of the motor at regular intervals and calculating a motor demand signal according to the control algorithm.

The motors used in the experimental apparatus are 3 phase brushless AC servo motors with electronic commutation and resolver feedback. By varying the voltage across the motor coils, the servo amplifiers induce a torque in the rotor of the motor by means of electromagnetic induction.

Previous work by Bradshaw (4,5) has claimed that the use of a PID controller with a velocity feedforward component based on an inverse model provides the best response for hybrid systems to date. Whilst this type of controller is not as robust as H∞ controllers, the predictable nature of the system implies that the feedforward component will produce low steady state errors and system stability.

### 3.1 Controller Tuning

The tuning of a single motor is a relatively simple process when compared to the tuning of a multi-axis machine, particularly where the inputs are mechanically coupled. In this case, the output from one motor effects the position of other motors in the machine and incorrect tuning can lead to ineffective machines with poor stability. Tuning a multi-axis machine, such as the hybrid five, is often done during the normal running operation. This is because the system can often be damaged by step inputs. The control algorithm used is of the following form;

$$V_{out} \propto K_P e_i + K_I \sum e_i + K_D (e_i - e_{i-1}) - K_V (p_i - p_{i-1}) + K_F (d_i - d_{i-1}) \qquad ...(5)$$

The actual scaling between error and output voltage, for proportional gain only, is shown below. The other gain terms have similar scaling factors.

$$V_{out} \propto e_i \times \frac{K_P}{256} \times \frac{10}{2048} \qquad ...(6)$$

### 3.2 Experimental Results

This section presents results for the purpose of evaluating the performance of the hybrid machine. The experiments consist of running the mechanism through a number of cycles and logging position and velocity data for both the CV and servo motors as well as observing the trajectory of the end effector with the LED and camera combination.

Figure 3.3 shows the actual end effector motion of the real machine in comparison to the precision points used to define the desired motion. This plot was obtained from long exposure photographs, where the trajectory of the end effector was traced by an LED. The trace obtained was copied onto acetate, and then projected and traced onto graph paper so that the precision points could be plotted. The desired start position of the mechanism is shown, as is the direction of rotation.

In addition to this, the following images were obtained by the use of a video camera. The mechanism was filmed around the cycle and a number of images were acquired. The purpose of these images is to show the actual end effector position for a given CV motor input angle and so compare this to the desired position. This differs from the information shown in Figure3.3 where the timing of the input is not shown. A sample image is given in Figure 3.4.

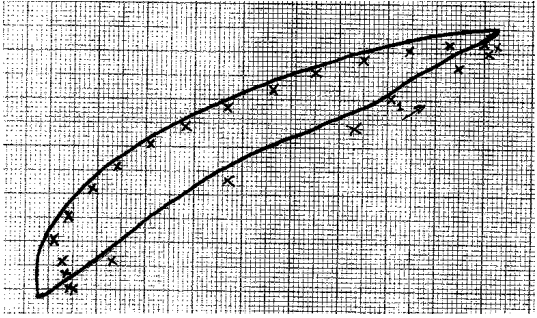
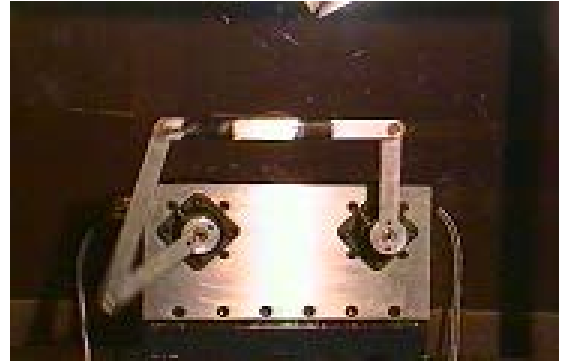

Figure 3.3 : End effector motion     Figure 3.4 : Machine position (217°)

Table II compares the actual co-ordinates of the end effector with the desired co-ordinates for three positions. The co-ordinates are expressed in millimetres relative to the CV motor position.

| CV INPUT | ACTUAL CO-ORDS (mm) | DESIRED CO-ORDS (mm) |
|---|---|---|
| 217° | -52,147 | -46,148 |
| 102° | 194,262 | 186,270 |
| 4° | 242,246 | 230,247 |

Table II : Comparison of end effector positions

The graph in Figure 3.5 shows the data logged from the machine for demand position, actual position and position error for the servo axis in units of resolver counts. The servo amplifiers are set to a resolution of 4096 counts per cycle (360°). By examining this data it can be seen that for most of the machine cycle, the actual position lags the demand position. This explains the relative position of the actual to the desired coupler curve observed in Figure 3.3.

The maximum position error occurs at the start of the return period of the motion. The value at this point is approximately -115 resolver counts, or approximately 10°. The maximum actual position of the servo motor is 1135 resolver counts. This is equivalent to 99.76°. This compares to the theoretical prediction of 99.78°. The demand and actual velocity of the servo axis, in units of resolver counts per second, were also used to assess the quality of the motion generation. There is a certain amount of "noise" on the demand signal, which is caused by the controller trying to correct position error. This causes the actual velocity signal to be quite oscillatory in nature.

The maximum actual velocity is 3072 resolver counts per second. This is equivalent to 4.71 radians per second. This compares to the theoretical prediction of 4.48 radians per second. Similarly, the minimum actual velocity is -12288 resolver counts per second. This is equivalent to -18.85 radians per second. This compares to the theoretical prediction of -18.96 radians per second.

Due to the uniform motion on the CV axis, the differences between demand and actual position are less significant. The graph in Figure 3.6 shows the position error only for this axis. The maximum value occurs at the beginning of the cycle and has a value of approximately 19 resolver counts, or 1.7°.

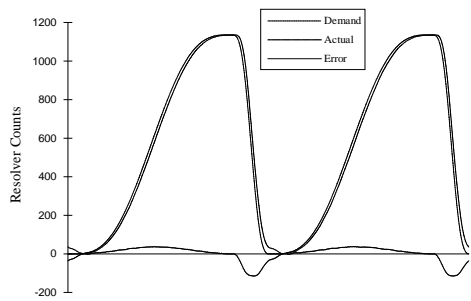 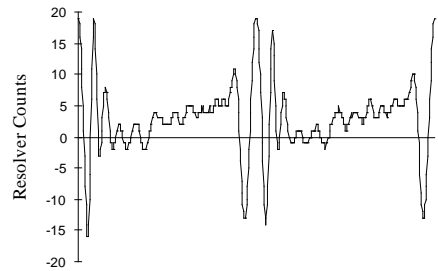

Figure 3.5 : Servo motor position data     Figure 3.6 : Position error on CV axis

The constant demand and the actual velocity for the CV axis were also considered. Again, the units are in resolver counts per second, where the demand velocity of 4096 cps is equal to 60 rpm. The maximum velocity error has a value of 1024 cps.

In addition to logging kinematic data concerning the mechanism to validate that the output motions are similar to the desired motions some simple power measurements were carried out to show the power distribution between CV and servo motor. Readings of RMS power were taken from a single phase of each motor supply, and averaged for several cycles of the mechanism. A three phase compensation was included which calculated the total power by referencing this phase to the remaining phases. The power drawn per cycle by the CV motor is 3.3429 kW, whilst the power drawn by the servo motor is 1.5696 kW. These figures can only be viewed as being approximations to the motor power consumption as they are actually measurements of the power drawn by the servo amplifier on each axis. However, they do indicate that the CV axis is contributing more to the overall motion than the servo axis.

## 4       DISCUSSION OF RESULTS

The results presented in this paper illustrate that the developed design methodology does produce results that can be turned into practical machines. However, several factors need to be highlighted. The most important of these is that the tuning of the controller is critical to producing the desired motions. The experimental tuning process used was halted as soon as the end effector of the mechanism appeared to generate the desired output motion, and there was no observable oscillation on either axis as the machine completed one cycle. The process was halted due to the tedious nature of finding optimal gain settings. The search for optimal gain settings is complex due to the effects that each gain term has on the output motion.

If an experimental tuning approach is to be used, then it is recommended that a more rigorous understanding of the effects of changing gains terms is obtained through careful experimentation. This would allow the controller to be tuned using tighter termination criteria which should produce better results, though the time taken to tune the controller for excellent performance may be excessive. Considering that the controller must be retuned if the servo motor motion profile is changed, then an experimental tuning approach is not recommended. However, by examining the graphs in Figure 3.6 and Figure 3.7 it can be seen that the position error is essentially cyclic. This implies that better control could be achieved by developing models of the hybrid machine. These models could be used to replace the current velocity feedforward component in the controller with an inverse model feedforward component as used by Bradshaw (4,5).

## 5  CONCLUSIONS

This paper has described the experimental apparatus developed to test that the developed design methodology produces solutions that can be turned into practical machines. Results have been presented which show that the real machine produces an adequate approximation to the desired motion. Errors between the desired and actual motion can be attributed to a set of non-optimal controller gains and the effects of a cross coupled system. The real machine shows a power distribution between CV and servo motors of approximately 2:1. This indicates that the bulk of the motion is being provided by the CV motor, whilst the servo motor is providing the flexibility to ensure that the desired coupler curve is being produced. This proves that the machine is acting as a hybrid device

## 6  REFERENCES

(1) **Tokuz, L.C.** *Hybrid Machine Modelling and Control*, Ph.D. Thesis, Liverpool Polytechnic, 1992
(2) **Tokuz, L.C. & Rees Jones, J.** *Programmable Modulation of Motion Using Hybrid Machines*, Institution of Mechanical Engineers, C414/071, pp 85-91, 1991
(3) **Greenough, J.D.** *Design of High Speed Machines*, Ph.D. Thesis, Liverpool John Moores University, 1996
(4) **Bradshaw, W.K.** *Control of Hybrid Machines*, Ph.D. Thesis, Liverpool John Moores University, 1997
(5) **Bradshaw, W.K, Gilmartin, M.J. & Douglas, S.S.** *Comparison of Three Algorithms for the Control of a Servo Motor*, Proceedings of the Ninth World Congress on the Theory of Machines and Mechanisms, Vol 2, pp 1355-1359, 1995
(6) **Connor, A.M.** *The Synthesis of Hybrid Mechanisms Using Genetic Algorithms*, Ph.D. Thesis, Liverpool John Moores University, 1996
(7) **Connor, A.M., Douglas, S.S. & Gilmartin, M.J.** *The Use of Harmonic Information in the Optimal Synthesis of Mechanisms*, Journal of Engineering Design, 9(3), 239-249
(8) **Connor, A.M., Douglas, S.S. & Gilmartin, M.J.** *The Synthesis of Path Generating Mechanisms Using Genetic Algorithms*, Proceedings of the 10th International Conference of Applications of Artificial Intelligence in Engineering, pp 237-244, 1995